\documentclass[11pt]{article}
\usepackage{amssymb}
\usepackage[top=23mm,bottom=23mm,left=23mm,right=23mm]{geometry}
\usepackage{graphicx}
\usepackage{amsmath}
\title{Transition density estimates for subordinated reflected Brownian motion on simple nested fractals}
\author{Hubert Balsam}

\usepackage{color}

\makeatletter
\@addtoreset{equation}{section}

\makeatother

\begin{document}

\newcommand{\rd}{\mathbb{R}^d}
\newcommand{\dist}{\emph{dist}}
\newcommand{\tr}{\mbox{tr}}
\newtheorem{theo}{\bf Theorem}[section]
\newtheorem{coro}{\bf Corollary}[section]
\newtheorem{lem}{\bf Lemma}[section]
\newtheorem{rem}{\bf Remark}[section]
\newtheorem{defi}{\bf Definition}[section]
\newtheorem{exam}{\bf Example}[section]
\newtheorem{fact}{\bf Fact}[section]
\newtheorem{prop}{\bf Proposition}[section]
\newtheorem{oq}{\bf Open question}

\maketitle

\begin{abstract}
In this paper we prove matching upper and lower  bounds for the transition density function of the subordinate reflected Brownian motion on  fractals.
\end{abstract}

\footnotetext{H. Balsam,  Department of Mathematics, Computer Science and Mechanics, University of Warsaw, ul. Banacha 2, 02-097 Warszawa, Poland \\
		{\it  e-mail: h.balsam@mimuw.edu.pl}}

\section{Introduction}

Stochastic processes on fractals, more generally on irregular sets, have been studied for over 40 years. The Brownian motion is the first process constructed in various spaces, such as the Sierpi\'{n}ski carpet \cite{bib:Bar-Bas}, the Sierpi\'{n}ski gasket \cite{bib:Bar-Per}, post critically-finite sets \cite{bib:Kum}, as well as on more general sets \cite{bib:Stu, bib:Lin,bib:Kum,bib:Fuk1}.

Suposse that the Brownian motion on an unbounded nested fractal $\mathcal{K}^{\left\langle \infty \right\rangle}$ has been constructed. For $M \in \mathbb{Z_+}$, we  construct the reflected Brownian motion on those compact fractals $\mathcal{K}^{\left\langle M \right\rangle}$  that obey the {\em good labeling property} (see Section \ref{defi:MC} for the definitions of $\mathcal{K}^{\left\langle \infty \right\rangle}, \mathcal{K}^{\left\langle M \right\rangle}$ and the \textit{good labeling property}). This has been achieved in \cite{bib:KOP} via a suitable projection procedure. The reflected Brownian motion on $\mathcal K^{\langle M\rangle}$  is a conservative diffusion process, whose transition density function satisfies  (see \cite{Olszewski}):
\begin{eqnarray*}
c_1t^{-\frac{d}{d_{w}}} \cdot {\rm e}^{-c_2\left(\frac{|x-y|}{t^{1/d_w}}\right)^\frac{d_w}{{d_{J}-1}}}  \leq & g_{M}(t,x,y) \leq & c_3 t^{-\frac{d}{d_{w}}} \cdot {\rm e}^{-c_4\left(\frac{|x-y|}{t^{1/d_w}}\right)^\frac{d_w}{{d_{J}-1}}} \mbox{ if } t<L^{Md_{w}}, \ x,y \in \mathcal{K}^{\left\langle M \right\rangle}\\[3mm]
c_5L^{-Md} \leq & g_{M}(t,x,y) \leq & c_6L^{-Md} \mbox{ if } t \geq L^{Md_{w}}, \ x,y \in \mathcal{K}^{\left\langle M \right\rangle},
\end{eqnarray*}
where $L$ is the scaling factor of $\mathcal{K}^{\left\langle \infty \right\rangle}$ and $\mathcal{K}^{\left\langle M \right\rangle},$ parameters $d, d_w, d_J$ depend on the geometry of the fractal, $c_1,\ldots c_6 > 0$ are absolute constants.
It is worth noting that when the time is large, then this transition density is comparable with $L^{-Md}$, which means that the process is roughly uniformly distributed over the $M-$ complex.

\

In this paper we would like to obtain estimates on the transition density function \linebreak for the subordinate reflected Brownian motion on $\mathcal{K}^{\left\langle M \right\rangle}.$ We will consider two classes of subordinate processes: $\alpha-$stable processes and  $\alpha-$stable relativistic processes. When $\mathcal{K}^{\left\langle M \right\rangle}$ is the Sierpi\'nski gasket, it has been proven in \cite{bib:KaPP2} that the subordination and the reflection commute, and this property holds in present case too. Therefore it convenient to understand the subordinate reflected process as the process subordinate to the reflected Brownian motion via the given subordinator (stable or relativistic). \linebreak In this paper, we prove the following result. For the transition density of the $\alpha-$stable reflected Brownian motion on $\mathcal{K}^{\left\langle M \right\rangle}$, denoted $p^M_S(t,x,y)$ (Theorem \ref{th:g_m_stable}): there exist  positive constants $B_1,B_2,B_3,B_4,$ such that for $M \in \mathbb{Z}_+$ and $x,y \in \mathcal{K}^{\left\langle M \right\rangle}$
\begin{eqnarray*}
B_1p_{S}(t,x,y)  \leq  p^M_S(t,x,y) \leq & B_2p_{S}(t,x,y) & \rm{ if } \textrm{ $t<L^{\alpha Md_{w}}$} \\[2mm]
B_3L^{-Md}  \leq  p^M_S(t,x,y)  \leq &  B_4L^{-Md} & \rm{ if } \textrm{ $t \geq L^{\alpha Md_{w}}.$} \\
\end{eqnarray*}
In the case of subordination via  the relativistic subordinator i.e. the relativistic $\alpha-$stable Brownian motion, we get that  for any $M \in \mathbb{Z}_+$, the transition density of the reflected relativistic $\alpha-$stable process, denoted $p_R^M(t,x,y)$ on $\mathcal K^{\langle M\rangle},$  satisfies (Theorem \ref{th:g_m_rel}):
\begin{enumerate}
\item[1)] for $ t \geq L^{Md_w},  \ x,y \in \mathcal{K}^{\left\langle M \right\rangle}$ there exist constants $C_1,C_2$ such that
\begin{equation}
C_1L^{-Md} \leq p^M_R(t,x,y) \leq C_2L^{-Md}
\end{equation}
\item[2)] for $t < L^{Md_w},   \ x,y \in \mathcal{K}^{\left\langle M \right\rangle}$ there exists constant $H_1$ such that
$$
p_R(t,x,y) \leq p_R^M(t,x,y) \leq p_R(t,H_1x,H_1y).
$$
\end{enumerate}
Taking into account the estimates on the relativistic stable process on $\mathcal{K}^{\left\langle \infty \right\rangle}$ from \cite{bib:Bal-Kpp} (cf. formulas (\ref{eq:p_r_t_geq_1}),(\ref{eq:p_r_t_leq_1}),(\ref{eq:p_r_t_geq_1_dxy_leq_1}) below)  we get that there exist constants $C_3,\ldots C_{12} > 0$ such that
\begin{enumerate}
\item[1)] for $1 \leq t < L^{Md_w},  \ x,y \in \mathcal{K}^{\left\langle M \right\rangle}$
\begin{multline}
C_3t^{-d/d_w} \exp\left\{-C_4\min\left(|x-y|^{\frac{d_{w}}{d_{J}}},\left({|x-y|}{t^{-\frac{1}{d_w}}}\right)^{\frac{d_{w}}{d_{J}-1}}\right)\right\} \leq p^M_R(t,x,y)\\
 \leq C_5t^{-d/d_w} \exp\left\{-C_6\min\left(|x-y|^{\frac{d_{w}}{d_{J}}},\left({|x-y|}{t^{-\frac{1}{d_w}}}\right)^{\frac{d_{w}}{d_{J}-1}}\right)\right\}
\end{multline}
\item[2)] for $t \in (0,1), |x-y| \geq 1$
\begin{equation}
C_7t{\rm e}^{-C_8|x-y|^{\frac{d_{w}}{d_{J}}}} \leq p^M_R(t,x,y) \leq C_9t{\rm e}^{-C_{10}|x-y|^{\frac{d_{w}}{d_{J}}}}
\end{equation}
\item[3)] for $t \in (0,1), |x-y| < 1$
\begin{equation}
C_{11}t^{\frac{-d}{\alpha d_{w}}}\left(\left(\frac{t^{\frac{1}{\alpha d_{w}}}}{|x-y|}\right)^{d +\alpha d_{w}} \wedge 1\right) \leq p^M_R(t,x,y) \leq C_{12}
 t^{\frac{-d}{\alpha d_{w}}}\left(\left(\frac{t^{\frac{1}{\alpha d_{w}}}}{|x-y|}\right)^{d +\alpha d_{w}} \wedge 1\right).
\end{equation}
\end{enumerate}

And again, these results show that those processes initially behave similarly to the 'original' ones and in large times they are almost uniformly distributed over entire $M-$complex.

\

The paper is organized as follows.  In Section 2 we provide definitions and notations regarding  unbounded simple nested fractals, subordination and reflected Brownian motion. Section 3 contains the proof of the estimates of the transition density for subordinated reflected Brownian motion via the  $\alpha-$stable subordinator, and Section 4 - via the relativistic subordinator.

\section{Preliminaries}
{\em Notation.}  Throughout the paper,  upper- and lowercase, numbered constants, $A_i, K_i, C_{i},c_{i}$ denote constants whose values, once fixed, will not change. Constants that are not numbered, i.e. $c, C, c', C',\ldots$  can change their value inside the proofs.
For two functions defined on a common domain,
$f\asymp g$ means that there is an absolute(independent of $t,x,y,M$) constant $C>0$ s.t. $\frac{1}{C} f(\cdot)\leq g(\cdot )\leq Cf(\cdot)$ , also $f \gtrsim g$ means that there is an absolute constant $C>0$ s.t. $ f(\cdot) \geq Cg(\cdot ).$

\subsection{Unbounded simple nested fractals} \label{sec:usnf}

The introductory part of this section follows the exposition of \cite{bib:Lin,bib:kpp-sausage,bib:kpp-sto}. Consider a collection of similitudes $\Psi_i : \mathbb{R}^2 \to \mathbb{R}^2$ with a common scaling factor $L>1,$  and a common isometry part $U,$ i.e. $\Psi_i(x) = (1/L) U(x) + \nu_i,$  where  $\nu_i \in \mathbb{R}^2$, $i \in \{1, ..., N\}.$ We shall assume $\nu_1 = 0$. Then
there exists a unique nonempty compact set $\mathcal{K}^{\left\langle 0\right\rangle}$ (called {\em the  fractal generated by the system} $(\Psi_i)_{i=1}^N$) such that $\mathcal{K}^{\left\langle 0\right\rangle} = \bigcup_{i=1}^{N} \Psi_i\left(\mathcal{K}^{\left\langle 0\right\rangle}\right)$.  As $L>1$, each similitude has exactly one fixed point and there are exactly $N$ fixed points of the transformations $\Psi_1, ..., \Psi_N$. Let $F$ be the collection of those fixed points.

\begin{defi}[\textbf{Essential fixed points}]
A fixed point $x \in F$ is an essential fixed point if there exists another fixed point $y \in F$ and two different similitudes $\Psi_i$, $\Psi_j$ such that $\Psi_i(x)=\Psi_j(y)$.
The set of all essential fixed points for transformations $\Psi_1, ..., \Psi_N$ is denoted by $V_{0}^{\left\langle 0\right\rangle}$, let $K=\# V^{\left\langle 0\right\rangle}_{0}$.
\end{defi}

\begin{figure}[ht]
\centering
	\includegraphics[scale=1]{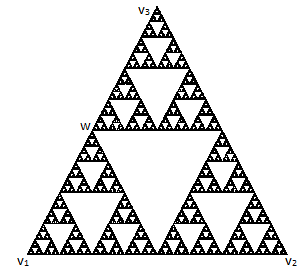}
\caption{Essential fixed points of the Sierpi\'nski triangle.}
\label{fig:essfix}
\end{figure}

\begin{exam}
The Sierpi\'nski triangle (Figure \ref{fig:essfix}) is constructed by 3 similitudes
$$
\Psi_1(x,y) =(\frac{x}{2},\frac{y}{2}), \mbox{ } \Psi_2(x,y) =(\frac{x}{2},\frac{y}{2}) + (\frac{1}{2},\frac{1}{2}), \mbox{ } \Psi_3(x,y) =(\frac{x}{2},\frac{y}{2}) + (\frac{1}{4},\frac{\sqrt{3}}{4})
$$
with scale factor $L = 2.$ The fixed points $v_i$ of the $\Psi_i'$s for $i=1,2,3$ are essential fixed points. For example, the vertex $v_1$ is an essential fixed point, because $\Psi_3(v_1) = \Psi_1(v_3)=w$.
\end{exam}

\begin{defi}[\textbf{Simple nested fractal}]
\label{def:snf}
 The fractal $\mathcal{K}^{\left\langle 0 \right\rangle}$ generated by the system $(\Psi_i)_{i=1}^N$ is called a \emph{simple nested fractal (SNF)} if the following conditions are met.
\begin{enumerate}
\item $\# V_{0}^{\left\langle 0\right\rangle} \geq 2.$
\item (Open Set Condition) There exists an open set $U \subset \mathbb{R}^2$ such that for $i\neq j$ one has\linebreak $\Psi_i (U) \cap \Psi_j (U)= \emptyset$ and $\bigcup_{i=1}^N \Psi_i (U) \subseteq U$.
\item (Nesting) $\Psi_i\left(\mathcal{K}^{\left\langle 0 \right\rangle}\right) \cap \Psi_j \left(\mathcal{K}^{\left\langle 0 \right\rangle}\right) = \Psi_i \left(V_{0}^{\left\langle 0\right\rangle}\right) \cap \Psi_j \left(V_{0}^{\left\langle 0\right\rangle}\right)$ for $i \neq j$.
\item (Symmetry) For $x,y \in V_{0}^{\left\langle 0\right\rangle},$ let $S_{x,y}$ denote the symmetry with respect to the line bisecting the segment $\left[x,y\right]$. Then
\begin{equation}
\forall i \in \{1,...,M\} \ \forall x,y \in V_{0}^{\left\langle 0\right\rangle} \ \exists j \in \{1,...,M\} \ S_{x,y} \left( \Psi_i \left(V_{0}^{\left\langle 0\right\rangle} \right) \right) = \Psi_j \left(V_{0}^{\left\langle 0\right\rangle} \right).
\end{equation}
\item (Connectivity) On the set $V_{-1}^{\left\langle 0\right\rangle}:= \bigcup_i \Psi_i \left(V_{0}^{\left\langle 0\right\rangle}\right)$ we define the graph structure $E_{-1}$ as follows:\\
$(x,y) \in E_{-1}$ if and only if $x, y \in \Psi_i\left(\mathcal{K}^{\left\langle 0 \right\rangle}\right)$ for some $i$.\\
Then the graph $(V_{-1}^{\left\langle 0\right\rangle},E_{-1} )$ is required to be connected.
\end{enumerate}
\end{defi}

If  $\mathcal{K}^{\left\langle 0 \right\rangle}$ is a simple nested fractal, then we let
\begin{equation} \label{eq:Kn}
\mathcal{K}^{\left\langle M\right\rangle} = L^M \mathcal{K}^{\left\langle 0\right\rangle}, \quad M \in \mathbb{Z},
\end{equation}
and
\begin{equation} \label{eq:Kinfty}
\mathcal{K}^{\left\langle \infty \right\rangle} = \bigcup_{M=0}^{\infty} \mathcal{K}^{\left\langle M\right\rangle}.
\end{equation}
The set $\mathcal{K}^{\left\langle \infty \right\rangle}$ is the \textbf{unbounded simple nested fractal (USNF)} we shall be working with (see \cite{bib:kpp-sausage}).  Its fractal (Hausdorff) dimension is equal to  $d=\frac{\log N}{\log L}$. The Hausdorff  measure in dimension $d$ will be denoted by $\mu$. It will be normalized to have $\mu\left(\mathcal{K}^{\left\langle 0\right\rangle}\right)=1$.

The remaining notions are collected in a single definition.
\begin{defi}\label{defi:MC}
Let $M\in\mathbb Z.$
\begin{itemize}
\item[(1)] $M$-complex: \label{def:Mcomplex}
every set $\Delta \subset \mathcal{K}^{\left\langle \infty \right\rangle}$ of the form
\begin{equation} \label{eq:Mcompl}
\Delta  = \mathcal{K}^{\left\langle M \right\rangle} + \nu_{\Delta},
\end{equation}
where $\nu_{\Delta}=\sum_{j=M+1}^{J} L^{j} \nu_{i_j},$ for some $J \geq M+1$, $\nu_{i_j} \in \left\{\nu_1,...,\nu_N\right\}$.
\item[(2)] Vertices of $\mathcal{K}^{\left\langle M \right\rangle}$:
$$
V^{\left\langle M\right\rangle}_{M} = V\left(\mathcal{K}^{\left\langle M \right\rangle}\right) = L^M V^{\left\langle 0\right\rangle}_{0}.
$$
\item[(3)] Vertices of all 0-complexes inside the unbounded nested fractal:
$$
V^{\left\langle \infty \right\rangle}_{0} = \bigcup_{M=0}^{\infty} V^{\left\langle M\right\rangle}_{0}.
$$
\item[(4)] Vertices of $M$-complexes from the unbounded fractal:
$$
V^{\left\langle \infty \right\rangle}_{M} = L^{M} V^{\left\langle \infty \right\rangle}_{0}
$$
\end{itemize}
\end{defi}

To define the reflected process, we need the good labeling property
introduced in \cite{bib:KaPP2}. We briefly sketch this idea.

\subsection{Good labelling and projections}\label{sec:labelling}
This section follows Section 3 of (\cite{bib:KaPP2}).\\

Recall that $K$ is the number of essential fixed points,  consider set of labels $\mathcal{A}:=\left\{a_1, a_2,a_3,...,a_K\right\}$ and a function  $l_M: V^{\left\langle \infty \right\rangle}_{M} \to \mathcal{A}.$ From \cite[Proposition 2.1]{bib:KaPP2} we have that there exist exactly $K$ different rotations $R_i$  around the barycenter of $\mathcal K^{\langle M \rangle}$,  mapping $V(\mathcal{K}^{\left\langle M \right\rangle})$ onto $V(\mathcal{K}^{\left\langle M \right\rangle}).$ Let us denote them as $\{R_1, ..., R_K\}=: \mathcal{R}_M.$\\
\begin{defi}\label{def:glp}
We say fractal that the fractal  $\mathcal K^{\langle \infty\rangle}$ has Good Labeling Property(GLP) if for some $M \in \mathbb{Z}$ there exist a function  $\ell_M: V^{\left\langle \infty \right\rangle}_{M} \to \mathcal{A}$ such that:
\begin{itemize}
\item[(1)] The restriction of $\ell_M$ to $V^{\left\langle M \right\rangle}_{M}$ is a bijection onto $\mathcal{A}$.
\item[(2)] For every $M$-complex $\Delta$ represented as
$$
\Delta  = \mathcal{K}^{\left\langle M \right\rangle} +\nu_{\Delta},
$$
where $\nu_{\Delta}=  \sum_{j=M+1}^{J} L^{j} \nu_{i_j},$  with some $J \geq M+1$ and $\nu_{i_j} \in \left\{\nu_1,...,\nu_N\right\}$ (cf. Def. \ref{def:Mcomplex}), there exists a rotation $R_{\Delta} \in \mathcal{R}_M$ such that
\begin{equation}\label{eq:glp-h}
\ell_M(v)=\ell_M\left(R_{\Delta}\left(v -\nu_{\Delta}\right)\right) , \quad v \in V\left(\Delta_M\right).
\end{equation}
\end{itemize}
\end{defi}
Now for the fractal $\mathcal K^{\langle \infty\rangle}$ having the GLP we define a projection map $\pi_{M}: \mathcal{K}^{\left\langle \infty \right\rangle} \rightarrow \mathcal{K}^{\left\langle M \right\rangle}$ as
$$
\pi_M(x) := R_{\Delta_M}\left(x -\nu_{\Delta_M}\right),\quad x\in \mathcal K^{\langle \infty\rangle},
$$
where $\Delta_M = \mathcal{K}^{\left\langle M \right\rangle} +  \sum_{j=M+1}^{J} L^{j} \nu_{i_j}$ is an $M$-complex containing $x$ represented as in definition (\ref{eq:glp-h}). More information regarding GLP and projections can be found in \cite{bib:KOP}.

\subsection{Stochastic processes on USNFs}

\subsubsection{Brownian motion on the unbounded fractal} \label{sec:free_BM}

Let $Z=(Z_t, \mathbf{P}^{x})_{t \geq 0, \, x \in \mathcal{K}^{\left\langle \infty \right\rangle}}$ be \emph{the Brownian motion} on the USNF $\mathcal{K}^{\left\langle \infty \right\rangle}$ \cite{bib:Kus2,bib:Lin}. It is a strong Markov, Feller process with transition probability densities $g(t,x,y)$ with respect to the $d$-dimensional Hausdorff measure $\mu$ on $\mathcal{K}^{\left\langle \infty \right\rangle}$, which are jointly continuous on $(0,\infty) \times \mathcal{K}^{\left\langle \infty \right\rangle} \times \mathcal{K}^{\left\langle \infty \right\rangle},$  satisfy the scaling property
$$
g(t,x,y) = L^{d_f} g(L^{d_w} t, L x, L y), \quad t>0, \ \ x, y \in \mathcal{K}^{\left\langle \infty \right\rangle}.
$$
and satisfy the subgaussian estimate
\begin{multline}
\label{eq:kum}
K_{1} t^{-d_s/2} \exp \left(-K_{2} \left(\frac{\left|x-y \right|^{d_w}}{t} \right)^{\frac{1}{d_J -1}} \right) \leq g(t,x,y) \\
\leq K_{3} t^{-d_s/2} \exp \left(-K_{4} \left(\frac{\left|x-y \right|^{d_w}}{t} \right)^{\frac{1}{d_J -1}} \right), \quad t>0, \ \ x,y \in \mathcal{K}^{\left\langle \infty \right\rangle},
\end{multline}
where $d_w$  is the walk dimension of $\mathcal{K}^{\left\langle \infty \right\rangle},$  $d_s=2d/d_w$ is its  spectral dimension,  and $d_J > 1$ is the so-called \emph{chemical exponent} of $\mathcal{K}^{\left\langle \infty \right\rangle}$. Constants $K_1, K_2, K_3, K_4$ are absolute. Typically $d_w \neq d_J$, but sometimes (e.g. for the Sierpi\'nski gasket) one has $d_w = d_J$, see \cite[Theorems 5.2, 5.5]{bib:Kum}.

\subsubsection{Reflected Brownian motions} \label{sec:reflected_BM}

Suppose now that the unbounded fractal $\mathcal K^{\langle\infty\rangle}$ has the GLP. For an arbitrary $M\in\mathbb Z$ the \emph{reflected Brownian motion} on $\mathcal{K}^{\left\langle M\right\rangle}$ is defined canonically by (see \cite{bib:Kus2})
\begin{equation}\label{eq:process-def}
Z_t^M = \pi_M(Z_t),
\end{equation}
where $\pi_M: \mathcal{K}^{\left\langle \infty  \right\rangle} \to \mathcal{K}^{\left\langle M\right\rangle}$ is the projection from in Section 3.\\
Its transition density $g_M(t,x,y):(0,\infty)\times \mathcal{K}^{\left\langle \infty \right\rangle} \times \mathcal{K}^{\left\langle M \right\rangle} \rightarrow (0,\infty) $ is given by
\begin{equation}\label{eq:g_m_defi}
g_M(t,x,y) \asymp \left\{ \begin{array}{ll}
\Sigma_{y' \in \pi^{-1}_M(y)}g(t,x,y')  & \mbox{ if } y \in \mathcal{K}^{\left\langle M \right\rangle} \char`\\ V_M^{\left\langle M \right\rangle} \\
& \\
\Sigma_{y' \in \pi^{-1}_M(y)} \cdot {\rm rank}(y')  & \mbox{ if } y \in V_M^{\left\langle M \right\rangle},\\
\end{array} \right.
\end{equation}
where ${\rm rank}(y_0)$ is the number of $M$-complexes meeting at the point $y_0 \in V_M$.\\
It has been proven in \cite{Olszewski} that the transition density of this process satisfies (Theorem  3.1 from \cite{Olszewski}):
\begin{align}\label{eq:olsz:g_m}
c_1(f_{c_2}(t,|x-y|) \vee h_{c_3}(t,M)) \leq g_{M}(t,x,y) \leq c_4(f_{c_5}(t,|x-y|) \vee h_{c_6}(t,M)),
\end{align}
where $c_1,c_2$ \ldots $c_6$ are certain nonnegative constants independent of $M$, and
\begin{eqnarray*}
f_{c}(t,r) & = & t^{-\frac{d}{d_{w}}} \cdot {\rm e}^{-c\left(\frac{r}{t^{1/d_w}}\right)^\frac{d_w}{{d_{J}-1}}},\\[2mm]
h_{c}(t,M) & = & L^{-dM}\left(\frac{L^{M}}{t^{1/d_{w}}} \vee 1\right) ^{d-\frac{d_{w}}{d_{J}-1}} \cdot {\rm e}^{ -c\left( \frac{L^{M}}{t^{1/d_{w}}} \vee 1\right)^{\frac{d_w}{{d_{J}-1}}}}.
\end{eqnarray*}
This estimate can be also written as (Corollary 3.1 of \cite{Olszewski}):
\begin{eqnarray}\label{eq:corollary_3.1_olsz}
c_1f_{c_2}(t,|x-y|)  \leq & g_{M}(t,x,y) \leq & c_3 f_{c_4}(t,|x-y|) \mbox{ if } t<L^{Md_{w}}, \ x,y \in \mathcal{K}^{\left\langle M \right\rangle}\\[3mm]
c_5L^{-Md} \leq & g_{M}(t,x,y) \leq & c_6L^{-Md} \quad \mbox{if} \quad t \geq L^{Md_{w}}, \ x,y \in \mathcal{K}^{\left\langle M \right\rangle}. \nonumber
\end{eqnarray}

\subsubsection{Subordinated reflected Brownian motion}
A subordinator $S=(S_{t},{\bf P})_{t \geq 0}$ is an  increasing L\'evy process on $[0, \infty)$ such that $S_{0}=0$ (see \cite{bib:Ber}). \linebreak The Laplace transform of its distribution  $\eta_{t}$ is given by
$$
 \int_{0}^{\infty}{\rm e}^{-\lambda s}\eta_{t}({\rm d}s)={\rm e}^{-t\phi(\lambda)}, \qquad  \lambda > 0.
$$

The function $\phi:(0,\infty) \rightarrow \mathbb{R},$ called \textit{the Laplace exponent of S,} can be expressed as\linebreak({L\'evy-Khintchine formula}):
$$
\phi(\lambda) = a\lambda + \int_{0}^{\infty}(1 - {\rm e}^{-\lambda x})\nu({\rm d}x)
$$
where $a\in\mathbb R$  is the {\em drift coefficient of $S,$ }and $\nu$ is the L\'evy measure of $S,$ i.e. a nonnegative, $\sigma-$finite,  Borel measure on $(0, \infty)$ such that:
\begin{equation}\label{eq:Levy_measure_cond}
\int_{0}^{\infty}(1 \land x)\nu({\rm d}x) < \infty.
\end{equation}

\noindent We will work with two classes of subordinators: $\alpha-$stable subordinators, with
\begin{equation}\label{eq:phi_stable}
\phi^S(\lambda) = \lambda^{\alpha} , \qquad \lambda > 0, \ \alpha \in (0,1)
\end{equation}
and relativistic $\alpha-$stable subordinators, with
\begin{equation}\label{eq:phi_relativistic}
\phi^R_m(\lambda) = (\lambda + m^{1/\alpha})^{\alpha} - m, \qquad \lambda,m > 0.
\end{equation}
 Denoting by $\eta_t(\cdot)$  the density of the $\alpha-$ stable subordinator and by $\eta_{t,m}(\cdot)$  the density of the relativistic $\alpha$-stable subordinator, we have  (see \cite[p. 3]{bib:Ryz}):
\begin{equation}\label{eq:density_relativistic}
 \eta_{t,m}(s):={\rm e}^{-m^{1/\alpha}s + mt}\eta_{t}(s), \qquad m > 0, \ \alpha \in (0,1), \quad s,t > 0.
\end{equation}

\subsubsection{Subordinated processes }

Assume that $(Z_{t},{\bf P}_{x})_{x \in \mathcal{K}^{\left\langle \infty \right\rangle}, t \geq 0}$ and $(Z^M_{t},{\bf P}_{x})_{x \in \mathcal{K}^{\left\langle M \right\rangle}, t \geq 0}$ is the Brownian motion on $\mathcal{K}^{\left\langle \infty \right\rangle}$ and $\mathcal{K}^{\left\langle M \right\rangle}$ respectively, and let $S$ be a subordinator independent of $Z$. We define the \emph{subordinate Brownian motion} $X = (X_t)_{t \geq 0}$ and the the \emph{subordinate reflected Brownian motion} $X^M = (X^M_t)_{t \geq 0}$ by
$$
X_t := Z_{S_t}, \quad t \geq 0,
$$
and
$$
X^M_t := Z^M_{S_t}, \quad t \geq 0
$$
respectively. These processes are c\`adl\`ag Markov processes with    transition densities given by:
\begin{equation}\label{eq:p_general_subordinated}
p(t,x,y) = \int_0^{\infty} g(u,x,y) \eta_t({\rm d}u), \quad  t>0, \ x,y, \in \mathcal{K}^{\left\langle \infty \right\rangle},
\end{equation}
and
\begin{equation}\label{eq:p_m}
p_M(t,x,y) = \int_0^{\infty} g_M(u,x,y) \eta_t({\rm d}u), \quad  t>0, \ x,y, \in \mathcal{K}^{\left\langle M \right\rangle}.
\end{equation}
Papers \cite{bib:Bal-Kpp} and  \cite{bib:BSS} were devoted to obtaining estimates for $\alpha-$stable and relativistic processes on $d-$sets. Nested fractals fall within this cathegory. More precisely, since $\mathcal{K}^{\left\langle \infty \right\rangle}$ is a $d-$set carrying a fractional diffusion $(Z_{t},{\bf P}_{x})_{x \in \mathcal{K}^{\left\langle \infty \right\rangle}, t \geq 0}$ and $(X_t)_{t \geq 0}$ is defined by subordination, then the following estimates hold true:

\begin{enumerate}
\item[(1)] For the $\alpha-$ stable process on $\mathcal{K}^{\left\langle \infty \right\rangle}$ (see \cite{bib:BSS}):
\begin{equation}\label{eq:stos-stable}
p_S(t,x,y)  \asymp  t^{\frac{-d}{\alpha d_{w}}}\left(\left(\frac{t^{\frac{1}{\alpha d_{w}}}}{|x-y|}\right)^{d +\alpha d_{w}} \wedge 1\right),  \ x,y \in \mathcal{K}^{\left\langle \infty \right\rangle}.
\end{equation}
\item[(2)] For the relativistic $\alpha-$ stable process on $\mathcal{K}^{\left\langle \infty \right\rangle}$ (see \cite{bib:Bal-Kpp}) there exist constants $A_1, \ldots A_{10} > 0$ such that:
\begin{itemize}
\item[(a)] for $t \geq 1,  \ x,y \in \mathcal{K}^{\left\langle \infty \right\rangle}$
\begin{multline}\label{eq:p_r_t_geq_1}
A_1t^{-d/d_w} \exp\left\{-A_2\min\left(|x-y|^{\frac{d_{w}}{d_{J}}},\left({|x-y|}{t^{-\frac{1}{d_w}}}\right)^{\frac{d_{w}}{d_{J}-1}}\right)\right\} \leq p_R(t,x,y)\\
 \leq A_3t^{-d/d_w} \exp\left\{-A_4\min\left(|x-y|^{\frac{d_{w}}{d_{J}}},\left({|x-y|}{t^{-\frac{1}{d_w}}}\right)^{\frac{d_{w}}{d_{J}-1}}\right)\right\}
\end{multline}
\item[(b)] for $t \in (0,1), |x-y| \geq 1$
\begin{equation}\label{eq:p_r_t_leq_1}
A_5t{\rm e}^{-A_6|x-y|^{\frac{d_{w}}{d_{J}}}} \leq p_R(t,x,y) \leq A_7t{\rm e}^{-A_8|x-y|^{\frac{d_{w}}{d_{J}}}}
\end{equation}
\item[(c)] for $t \in (0,1), |x-y| < 1$
\begin{equation}\label{eq:p_r_t_geq_1_dxy_leq_1}
A_9t^{\frac{-d}{\alpha d_{w}}}\left(\left(\frac{t^{\frac{1}{\alpha d_{w}}}}{|x-y|}\right)^{d +\alpha d_{w}} \wedge 1\right) \leq p_R(t,x,y) \leq A_{10}
 t^{\frac{-d}{\alpha d_{w}}}\left(\left(\frac{t^{\frac{1}{\alpha d_{w}}}}{|x-y|}\right)^{d +\alpha d_{w}} \wedge 1\right).
\end{equation}
\end{itemize}
\end{enumerate}
Therefore one already has estimates for stable and relativistic processes on the infinite fractal $\mathcal{K}^{\left\langle \infty \right\rangle}.$
Now we are ready to formulate and prove corresponding estimates for the reflected processes on $\mathcal{K}^{\left\langle M \right\rangle}$.
\noindent

\section{Transition density estimate for the reflected $\alpha$-stable Brownian motion on $\mathcal{K}^{\left\langle M \right\rangle}$.}
We stand with the simpler case of the reflected $\alpha-$stable Brownian motion, obtained by  the subordination of the reflected Brownian motion on $\mathcal{K}^{\left\langle \infty \right\rangle}.$  Let $M \in \mathbb{Z_+}$ be fixed. We have the following.
\begin{theo}\label{th:g_m_stable}
Let $X_{t}$ be the $\alpha$-stable reflected process on $\mathcal{K}^{\left\langle M \right\rangle}$, with density function $p^M_S(\cdot, \cdot, \cdot)$ given by  (\ref{eq:p_m}). Then there exist constants $B_1, B_2, B_3, B_4 > 0$ independent of $M$ such that
\begin{eqnarray*}
B_1p_{S}(t,x,y)  \leq  p^M_S(t,x,y) \leq & B_2p_{S}(t,x,y) & \rm{ if } \textrm{ $t<L^{\alpha Md_{w}}$} \\[2mm]
B_3L^{-Md}  \leq  p^M_S(t,x,y)  \leq & B_4L^{-Md} & \rm{ if } \textrm{ $t \geq L^{\alpha Md_{w}}.$} \\
\end{eqnarray*}
\end{theo}
\noindent{\bf Proof.}\\
Before we begin, observe that
for any $c > 0$ the function $f_c(t,r)=t^{-\frac{d}{d_{w}}} \cdot {\rm e}^{-c\left(\frac{r}{t^{1/d_w}}\right)^\frac{d_w}{{d_{J}-1}}}$ is monotone decreasing in $r$, therefore, since for $x,y \in \mathcal{K}^{\left\langle M \right\rangle}, |x-y| \leq L^{M}, $ we have for any $c > 0$
$$
f_{c}(t,|x-y|) \geq f_{c}(t,L^{M})
$$
and also for any given constants $c_1, c_2, c_3 >0, c_1 < c_2$ there exists constants $c(c_{1},c_{2},c_{3}),c'(c_{1},c_{2},c_{3})$ such that for all  $ s \in [ c_{1} L^{Md_{w}},c_{2}L^{Md_{w}}]$ and $x,y \in \mathcal{K}^{\left\langle M \right\rangle}$ we have:
\begin{equation}\label{eq:f estimate}
c(c_{1},c_{2},c_{3})L^{-Md} \leq f_{c_{3}}(s,|x-y|) \leq c'(c_{1},c_{2},c_{3})L^{-Md}.
\end{equation}
Indeed:
$$
f_{c_{3}}(s,|x-y|) \leq (c_{1} \cdot L^{Md_{w}})^{-\frac{d}{d_{w}}} = c_{1}^{-\frac{d}{d_{w}}}L^{-Md}
$$
and given that $x \mapsto x^{-d/d_w}$ is decreasing and $x \mapsto {\rm e}^{-c\left(\frac{r}{x^{1/d_w}}\right)^\frac{d_w}{{d_{J}-1}}} $ is increasing (for $x > 0$) we get:
$$
f_{c_{3}}(s,|x-y|) \geq (c_{2}L^{Md_w})^{\frac{-d}{d_{w}}}f_{c_{3}}(c_{1}L^{Md_{w}},L^{M}) = c_{2}^{\frac{-d}{d_{w}}}{\rm e}^{-c_{3}c_{1}^{\frac{-d_w}{d_J-1}}}L^{-Md}.
$$
We now pass to the actual estimate.$\medskip$\\
\noindent\textsc{CASE 1. $t \geq L^{\alpha Md_{w}}.$} Let
$$
B_5 = (\frac{c_4}{K_2})^{\frac{d_J-1}{d_w}}.
$$
Using (\ref{eq:corollary_3.1_olsz}) we have (recall that $c_3, c_4$ do not depend on $M$):
\begin{eqnarray*}\label{upper t ver st}
p^M_S(t,x,y) & = & \int_0^{\infty} g_M(s,x,y) \eta_t({\rm d}s)\\
 & \leq & c_3\int_{0}^{L^{Md_{w}}}f_{c_4}(s,|x-y|)\eta_{t}(s){\rm d}s + c_6\int_{L^{Md_{w}}}^{\infty}L^{-Md}\eta_{t}(s){\rm d}s \\
 & \leq & \frac{c_3}{c_1}\int_{0}^{\infty}g(s,B_5|x-y|)\eta_{t}(s){\rm d}s + c_6\int_{0}^{\infty}L^{-Md}\eta_{t}(s){\rm d}s \\
 & \leq & cp_{S}(t,B_5x,B_5y) + c_6L^{-Md}
\end{eqnarray*}
as $\int_{0}^{\infty}\eta_{t}(s){\rm d}s = 1.\bigskip$\\
Given $t \geq L^{\alpha Md_{w}}$ we get from (\ref{eq:stos-stable}):
$$
p_{S}(t,B_5x,B_5y) \leq c\cdot t^{\frac{-d}{\alpha d_{w}}}\cdot (B_5^{-1} \wedge 1) \leq c \cdot L^{-Md}
$$
which means that
$$
p^M_S(t,x,y) \leq c \cdot L^{-Md}.
$$
In \cite[formula (10), p.4]{bib:BSS} we have that:
$$
\eta_{t}(u) \geqslant ctu^{-1-\alpha} \qquad {\rm for} \ t >0, \ u > u_{0}t^{1/\alpha}. \mbox { Without loss of generality } u_0 \geq 1.
$$
Due to $t \geq L^{\alpha Md_{w}}$ we have  $u_{0}t^{1/\alpha} \geq u_{0}L^{Md_{w}} \geq L^{Md_{w}}$ so using the subordination and (\ref{eq:corollary_3.1_olsz}):
$$
p^M_S(t,x,y)  = \int_{0}^{\infty}g_M(u,x,y)\eta_{t}(u){\rm d}u \geq   c_3\int_{u_{0}t^{1/\alpha}}^{\infty}L^{-Md}\eta_{t}(u){\rm d}u \geq  cL^{-Md}\int_{u_{0}t^{1/\alpha}}^{\infty}tu^{-1-\alpha}{\rm d}u =cL^{-Md}
$$
which  implies that for $t \geq L^{\alpha Md_{w}}$ the proof is done.$\medskip$\\
\noindent\textsc{CASE 2. $t < L^{\alpha Md_{w}}.$} Firstly, let us note that for any constant $A > 0$ there exists a constant $C(A)$ such that
$$
\frac{1}{C(A)}p_{S}(t, x,y) \leq p_{S}(t, Ax,Ay) \leq C(A)p_{S}(t, x,y).
$$
Let
$$
B_6 = (\frac{c_2}{K_2})^{\frac{d_J-1}{d_w}}.
$$
From (\ref{eq:olsz:g_m}):
\begin{eqnarray}\label{eq:pm-ps} \nonumber
p^M_S(t,x,y) & = & \int_{0}^{\infty}g_M(s,x,y)\eta_{t}(s){\rm d}s\\ [2mm] \nonumber
& \geq & c_1\int_{0}^{\infty}f_{c_2}(s,|x-y|)\eta_{t}(s){\rm d}s\\ [2mm] \nonumber
& \geq & c\int_{0}^{\infty}g(s,B_6x,B_6y)\eta_{t}(s){\rm d}s \\[2mm] \nonumber
& = &cp_{S}(t, B_6x,B_6y)\\[2mm]
& \geq & C(B_6)p_{S}(t, x,y)
\end{eqnarray}
so it is enough to show that
$$
p^M_S(t,x,y) \leq cp_{S}(t, x,y), \mbox{ with some $c$  independent of $M.$}
$$
Using (\ref{eq:corollary_3.1_olsz}) we have:
\begin{eqnarray*}\label{upper t ver st}
p^M_S(t,x,y) & \leq & c_3\int_{0}^{L^{Md_{w}}}f_{c_4}(s,|x-y|)\eta_{t}(s){\rm d}s + c_6\int_{L^{Md_{w}}}^{\infty}L^{-Md}\eta_{t}(s){\rm d}s \\
 & \leq & c_3\int_{0}^{\infty}f_{c_4}(s,|x-y|)\eta_{t}(s){\rm d}s + c_6\int_{L^{Md_w}}^{\infty}L^{-Md}\eta_{t}(s){\rm d}s \\
\end{eqnarray*}
 so that
\begin{equation}\label{eq:ps - I1}
p_{M,S}(t,x,y) \leq cp_{S}(t,x,y) + c_6\int_{L^{Md_w}}^{\infty}L^{-Md}\eta_{t}(s){\rm d}s.
\end{equation}
In \cite[formula (9), p.4]{bib:BSS} we have that $\eta_{t}(u) \leq ctu^{-1-\alpha}, \ \mbox{ for } t,u >0,$ so:
$$
I_{1}:= c_6\int_{L^{Md_{w}}}^{\infty}L^{-Md}\eta_{t}(s){\rm d}s  \leq ctL^{-M(d+\alpha d_{w})} \leq cL^{-Md}
$$
but, if $t^{\frac{1}{\alpha d_{w}}} \geq |x-y|,$ then
$$
p_{S}(t,x,y) \geq ct^{\frac{-d}{\alpha d_{w}}} \geq cL^{-Md} \geq cI_{1}.
$$
On the other hand, if $t^{\frac{1}{\alpha d_{w}}} < |x-y|$ then (as $x,y \in \mathcal{K}^{\left\langle M \right\rangle}$ i.e $|x-y|  \leq L^{M})$:
$$
p_{S}(t,x,y) \geq c\cdot \frac{t}{|x-y|^{d+\alpha d_{w}}} \geq ctL^{-M(d+\alpha d_{w})} \geq cI_{1},
$$
which means that $I_1 \lesssim p_S(t,x,y).$ From (\ref{eq:ps - I1}) we get $p^M_S(t,x,y) \lesssim p_{S}(t,x,y)$ and returning to (\ref{eq:pm-ps}) we have $p_{M,S}(t,x,y) \asymp p_{S}(t, x,y).$ The proof is done.\hfill$\Box$

\section{Transition density estimate for the $\alpha$-stable relativistic reflected Brownian motion}
In this section we provide the estimate of the density transition  of the  reflected Brownian motion obtained via the relativistic subordinator, i.e. $p_M(t,x,y) = \int_0^{\infty} g_M(s,x,y) \eta_{t,m}({\rm d}s)$  where $\eta_{t,m}$ is given by (\ref{eq:density_relativistic}).
\begin{theo}\label{th:g_m_rel}
Let $X_{t}$ be the relativistic $\alpha$-stable reflected process on $\mathcal{K}^{\left\langle M \right\rangle}$, with density function $p^M_R(\cdot, \cdot, \cdot)$. Then
\begin{enumerate}
\item[1)] for $ t \geq L^{Md_w},  \ x,y \in \mathcal{K}^{\left\langle M \right\rangle}$ there exist constants $C_1,C_2 > 0$ such that
\begin{equation}
C_1L^{-Md} \leq p_R^M(t,x,y) \leq C_2L^{-Md}.
\end{equation}
\item[2)] for $t < L^{Md_w},$ there exists constant $H_2$ such that
$$
p_R(t,x,y) \leq p_R^M(t,x,y) \leq p_R(t,H_2x,H_2y).
$$
\end{enumerate}
\end{theo}
This statement means that there exist constants $C_3,\ldots C_{12} > 0$ such that
\begin{enumerate}
\item[1)] for $1 \leq t < L^{Md_w},  \ x,y \in \mathcal{K}^{\left\langle M \right\rangle}$
\begin{multline}
C_3t^{-d/d_w} \exp\left\{-C_4\min\left(|x-y|^{\frac{d_{w}}{d_{J}}},\left({|x-y|}{t^{-\frac{1}{d_w}}}\right)^{\frac{d_{w}}{d_{J}-1}}\right)\right\} \leq p_R^M(t,x,y)\\
 \leq C_5t^{-d/d_w} \exp\left\{-C_6\min\left(|x-y|^{\frac{d_{w}}{d_{J}}},\left({|x-y|}{t^{-\frac{1}{d_w}}}\right)^{\frac{d_{w}}{d_{J}-1}}\right)\right\}
\end{multline}
\item[2)] for $t \in (0,1), |x-y| \geq 1$
\begin{equation}
C_7t{\rm e}^{-C_8|x-y|^{\frac{d_{w}}{d_{J}}}} \leq p_R^M(t,x,y) \leq C_9t{\rm e}^{-C_{10}|x-y|^{\frac{d_{w}}{d_{J}}}}
\end{equation}
\item[3)] for $t \in (0,1), |x-y| < 1$
\begin{equation}
C_{11}t^{\frac{-d}{\alpha d_{w}}}\left(\left(\frac{t^{\frac{1}{\alpha d_{w}}}}{|x-y|}\right)^{d +\alpha d_{w}} \wedge 1\right) \leq p_R^M(t,x,y) \leq C_{12}
 t^{\frac{-d}{\alpha d_{w}}}\left(\left(\frac{t^{\frac{1}{\alpha d_{w}}}}{|x-y|}\right)^{d +\alpha d_{w}} \wedge 1\right).
\end{equation}
\end{enumerate}

\noindent{\bf Proof.}
Let
$$
H_1 = (\frac{c_4}{K_2})^{\frac{d_J-1}{d_w}}, \
H_2 := \left(\frac{c_4}{K_4}\right)^{\frac{d_J-1}{d_w}}.
$$
These constants will be needed later.\medskip\\
\textsc{CASE 1.} $t \geq L^{Md_{w}}.$
Using (\ref{eq:corollary_3.1_olsz}) we have:
\begin{eqnarray*}
p^M_R(t,x,y) & = & \int_0^{\infty} g_M(s,x,y) \eta_{t,m}(s){\rm d}s\\[2mm]
& \leq & c_3\int_{0}^{L^{Md_{w}}}f_{c_4}(s,|x-y|)\eta_{t,m}(s){\rm d}s + c_6\int_{L^{Md_{w}}}^{\infty}L^{-Md}\eta_{t,m}(s){\rm d}s\\[2mm]
 & \leq & c_3{\rm e}^{mt}\int_{0}^{\infty}f_{c_4}(s,|x-y|){\rm e}^{-m^{\frac{1}{\alpha}}s}\eta_{t}(s){\rm d}s + c_6{\rm e}^{mt}\int_{0}^{\infty}L^{-Md}{\rm e}^{-m^{\frac{1}{\alpha}}s}\eta_{t}(s){\rm d}s\\[2mm]
 & =: & c_3I_2 + c_6{\rm e}^{mt}L^{-Md}I_3.
\end{eqnarray*}
We have $I_2 \leq p_R(t, H_1x,H_1y)$ and since $t \geq L^{Md_{w}}\geq 1$ we get from (\ref{eq:p_r_t_geq_1})
$$
I_2 \leq cL^{-Md}.
$$
The integral $I_3$ is the Laplace transform of $\eta_t$ evaluated at the at the  point $m^{\frac{1}{\alpha}},$ so
$$
I_3 = {\rm e}^{-t\phi^S(m^{\frac{1}{\alpha}})}={\rm e}^{-mt}.
$$
Altogether, there is a universal constant $a_1$ such that
\begin{equation}\label{eq:upper t large}
p^M_R(t,x,y) \leq a_1L^{-Md}, \mbox{ for } t \geq L^{Md_{w}},\ x,y  \in  \mathcal{K}^{\left\langle M \right\rangle}.
\end{equation}
The upper bound is done.$\medskip$\\
Now the lower bound. In the course of the proof of Theorem 3.1 in \cite{bib:Bal-Kpp}, p.193, we have proven that for any $m > 0$ there exist constants $L_1 = L_1(m), L_2 = L_2(m) > 1, a = a(m)$ such that

$$
\int_{L_{1}t}^{L_{2}t}{\rm e}^{-m^{1/\alpha}s}\eta_{t}(s){\rm d} s \geq a{\rm e}^{-mt}, \mbox{ for }  t > 0.
$$
Clearly
\begin{eqnarray*}
p^M_R(t,x,y) & = & \int_0^{\infty} g_M(s,x,y) \eta_{t,m}(s){\rm d}s\\[2mm]
& \geq & {\rm e}^{mt}\int_{L_1t}^{L_2t}g_M(s,x,y){\rm e}^{-m^{\frac{1}{\alpha}}s}\eta_{t}(s){\rm d}s
\end{eqnarray*}

We have two posibilities: $(1) \ L^{Md_w} \leq L_1t; \ (2) \ L_1t < L^{Md_w} < L_2t.$\\
$\bullet$ If $L^{Md_{w}}\leq L_{1}t < L_{2}t$, then from (\ref{eq:corollary_3.1_olsz}) and (\ref{eq:p_r_t_geq_1}):
$$
p^M_R(t,x,y) \geq c_5L^{-Md}{\rm e}^{mt}\int_{L_{1}t}^{L_{2}t}{\rm e}^{-m^{\frac{1}{\alpha}}s}\eta_{t}(s){\rm d}s \geq c_5a(m)L^{-Md}.
$$
$\bullet$ If $ L_{1}t < L^{Md_{w}} \leq L_{2}t$ which, in light of the assumption $t \geq L^{Md_w}$ requires $L_1 < 1,$ then from (\ref{eq:corollary_3.1_olsz}) and (\ref{eq:p_r_t_geq_1}):
$$
p^M_R(t,x,y) \geq c {\rm e}^{mt}\left(\int_{L_{1}t}^{L^{Md_{w}}}f_{A_2}(s,|x-y|){\rm e}^{-m^{1/\alpha}s}\eta_{t}(s){\rm d} s + \int_{L^{Md_{w}}}^{L_{2}t}L^{-Md}{\rm e}^{-m^{1/\alpha}s}\eta_{t}(s){\rm d} s\right),
$$
but if $ t \geq L^{Md_{w}}$ and $ L_{1}t < L^{Md_{w}}$ i.e. $L_1t \in[L_1L^{Md_w},L^{Md_w}]$ then there exists $L_{3}(t) \in  [L_{1},1)$ such that $L_{1}t = L_{3}L^{Md_{w}}.$ From (\ref{eq:f estimate}) with $c_1 = L_3(t), c_2 = 1, c_3 = A_2$ we get:
$$
f_{A_2}(s,|x-y|) \geq {\rm e}^{-A_2L_3(t)^{\frac{-d_w}{d_J-1}}}L^{-Md} \geq {\rm e}^{-A_2L_1^{\frac{-d_w}{d_J-1}}}L^{-Md}
$$
and further
\begin{eqnarray*}
p^M_R(t,x,y) & \geq  & c{\rm e}^{mt}\left(\int_{L_{1}t}^{L^{Md_{w}}}L^{-Md}{\rm e}^{-m^{1/\alpha}s}\eta_{t}(s){\rm d} s + \int_{L^{Md_{w}}}^{L_{2}t}L^{-Md}{\rm e}^{-m^{1/\alpha}s}\eta_{t}(s){\rm d} s\right) \\[2mm]
& = & cL^{-Md}{\rm e}^{mt}\int_{L_{1}t}^{L_{2}t}{\rm e}^{-m^{\frac{1}{\alpha}}s}\eta_{t}(s){\rm d}s \\[2mm]
& \geq & c\cdot a(m) L^{-Md}
\end{eqnarray*}
So we have just shown that
$$
p^M_R(t,x,y) \geq  cL^{-Md}
$$
so that for $t \geq L^{Md_{w}}$ the proof is complete.\\

\noindent\textsc{CASE 2.} $t < L^{Md_{w}}.$ From (\ref{eq:g_m_defi}) we have:
\begin{eqnarray*}
p^M_R(t,x,y) & = & \int_0^{\infty} g_M(s,x,y) \eta_{t,m}(s){\rm d}s\\[2mm]
& \geq & \int_0^{\infty} \Sigma_{y' \in \pi^{-1}_M(y)}g(s,x,y') \eta_{t,m}(s){\rm d}s\\[2mm]
& \geq & \int_0^{\infty} g(s,x,y) \eta_{t,m}(s){\rm d}s = p_R(t,x,y)
\end{eqnarray*}
so it will be enough to show the upper bound if we show that it holds
$$
p^M_R(t,x,y) \leq c\cdot p_R(t, H_2x,H_2y).
$$
Let
$$
K := 2m^{-\frac{1}{\alpha}+1}.
$$
Observe that, given (\ref{eq:f estimate}) we can adjust the estimate in (\ref{eq:corollary_3.1_olsz}) in such a way that the threshold is $KL^{Md_w}$ ($K$ remains fixed) and this only requires changes in constants. For simplicity, assume that the constants in (\ref{eq:corollary_3.1_olsz}) work for this threshold.\\
Since for $s \geq KL^{Md_w}$, it holds that $g_M(s,x,y) \leq cL^{-Md}$, so we have from (\ref{eq:corollary_3.1_olsz}) and (\ref{eq:f estimate}):
\begin{eqnarray*}
p^M_R(t,x,y) & = & \int_0^{\infty} g_M(s,x,y) \eta_{t,m}(s){\rm d}s\\[2mm]
 & \leq & {\rm e}^{mt}\int_{0}^{K\cdot L^{Md_w}}g_M(s,x,y){\rm e}^{-m^{\frac{1}{\alpha}}s}\eta_{t}(s){\rm d}s + c{\rm e}^{mt}\int_{K\cdot L^{Md_w}}^{\infty}L^{-Md}{\rm e}^{-m^{\frac{1}{\alpha}}s}\eta_{t}(s){\rm d}s\\[2mm]
 & \leq & c_3\int_{0}^{\infty}f_{c_4}(s,|x-y|)\eta_{t,m}(s){\rm d}s + c{\rm e}^{mt}\int_{K\cdot L^{Md_w}}^{\infty}L^{-Md}{\rm e}^{-m^{\frac{1}{\alpha}}s}\eta_{t}(s){\rm d}s\\[2mm]
 & \leq & c_3{\rm e}^{mt}\int_{0}^{\infty}f_{c_4}(s,|x-y|){\rm e}^{-m^{\frac{1}{\alpha}}s}\eta_{t}(s){\rm d}s + c_6{\rm e}^{mt}\int_{K\cdot L^{Md_w}}^{\infty}L^{-Md}{\rm e}^{-m^{\frac{1}{\alpha}}s}\eta_{t}(s){\rm d}s\\[2mm]
 & \leq & cp_R(t,H_2x,H_2y) + c{\rm e}^{mt}\int_{K\cdot L^{Md_{w}}}^{\infty}L^{-Md}{\rm e}^{-m^{\frac{1}{\alpha}}s}\eta_{t}(s){\rm d}s.
\end{eqnarray*}
Let
$$
I_4 := {\rm e}^{mt}\int_{K\cdot L^{Md_{w}}}^{\infty}L^{-Md}{\rm e}^{-m^{\frac{1}{\alpha}}s}\eta_{t}(s){\rm d}s.
$$
Due to the fact $\eta_{t}(u) \leq ctu^{-1-\alpha}, t,u > 0$ \cite[formula (9), p.4]{bib:BSS} then:
\begin{eqnarray*}
I_4 & \leq & cL^{-Md}{\rm e}^{mt}\int_{K\cdot L^{Md_{w}}}^{\infty}{\rm e}^{-m^{\frac{1}{\alpha}}s}ts^{-1-\alpha}{\rm d}s \\
 & \leq & cL^{-Md}  \cdot \frac{t}{(L^{Md_{w}})^{\alpha + 1}} {\rm e}^{mt}\int_{2m^{-\frac{1}{\alpha}+1}\cdot L^{Md_{w}}}^{\infty}{\rm e}^{-m^{\frac{1}{\alpha}}s}{\rm d}s \\
 & = & cL^{-Md}  \cdot \frac{t}{(L^{Md_{w}})^{\alpha + 1}} {\rm e}^{mt}{\rm e}^{-2mL^{Md_{w}}} \\
 & \leq & cL^{-Md}  \cdot \frac{t}{(L^{Md_{w}})^{\alpha + 1}} {\rm e}^{-mL^{Md_{w}}}.
\end{eqnarray*}
As $t < L^{Md_{w}}$, we have $\frac{t}{(L^{Md_{w}})^{\alpha + 1}} \leq 1$, so that
$$
I_4 \leq cL^{-Md}e^{-mL^{Md_w}},
$$
and further since $|x-y| \leq L^{M}$ and $d_J > 1$, we conclude with:\\
$\bullet$ for $ t\geq 1$
\begin{eqnarray*}
p_R(t,H_2x,H_2y) & \geq & ct^{-\frac{d}{d_{w}}}{\rm e}^{-K_6\left(\frac{c_4}{K_4}\right)^{\frac{d_J-1}{d_J}}|x-y|^{\frac{d_{w}}{d_{J}}}}\\
 & \geq & cL^{-Md}{\rm e}^{-K_6\left(\frac{c_4}{K_4}\right)^{\frac{d_J-1}{d_J}}(L^{M})^{\frac{d_{w}}{d_{J}}}}  \geq  cI_4
\end{eqnarray*}
for  $M \geq M_0$ where $M_0 = M_0(c_4, K_4, K_6, m, d_w, d_J)$ is an integer.\\
$\bullet$ for $ t \in (0,1), H_2|x-y| \geq 1$, we have:
\begin{eqnarray*}
p_R(t,H_2x,H_2y) & \geq & cte^{-K_{10}\left(\frac{c_4}{K_4}\right)^{\frac{d_J-1}{d_J}}|x-y|^{\frac{d_{w}}{d_{J}}}}\\
 & \geq & ct\cdot 1 \cdot {\rm e}^{-K_{10}\left(\frac{c_4}{K_4}\right)^{\frac{d_J-1}{d_J}}(L^{M})^{\frac{d_{w}}{d_{J}}}}  \geq  cI_4
\end{eqnarray*}
again, for $M \geq M_0'$ where $M_0' := M_0'(c_4, K_4, K_{10}, m, d_w, d_J)$ is an integer.\\
$\bullet$ for $ t \in (0,1), H_2|x-y| < 1$, we have:
\begin{eqnarray*}
p_R(t,H_2x,H_2y) & \geq & ct((H_1|x-y|)^{-d -\alpha d_{w}} \wedge t^{\frac{-d-\alpha d_{w}}{\alpha d_{w}}})\\
 & \geq & ct\cdot 1  \geq  cI_4
\end{eqnarray*}
which completes the proof of the theorem.\hfill$\Box$

\begin{center}
\section*{Acknowledgements}
\end{center}
Thank you to my supervisor Katarzyna Pietruska-Pa\l uba  for providing guidance throughout this paper and many valuable remarks.

\end{document}